\documentclass[12pt]{amsart}
\usepackage[all]{xy}
\usepackage{amssymb}

\newcommand{\CC}{\mathbb C}
\newcommand{\ZZ}{\mathbb Z}
\newcommand{\QQ}{\mathbb Q}
\newcommand{\PP}{\mathbb P}
\newcommand{\GG}{\mathfrak G}
\newcommand{\KZ}{KZ}
\newcommand{\llambda}{\lambda}
\newcommand{\Bun}{\mathcal Bun}
\newcommand{\Conf}{\mathcal Conf}
\newcommand{\OO}{\mathcal O}
\newcommand{\D}{\mathcal D}
\newcommand{\F}{\mathcal F}
\renewcommand{\P}{\mathcal P}
\newcommand{\LL}{\mathcal L}
\newcommand{\ad}{\mathop{\rm ad}\nolimits}
\newcommand{\rk}{\mathop{\rm rk}\nolimits}
\newcommand{\asym}{\mathop{\rm asym}\nolimits}
\newcommand{\Pic}{\mathop{\rm Pic}\nolimits}
\newcommand{\Ext}{\mathop{\rm Ext}\nolimits}

\begin{document}

\title
{Quantum Langlands duality and conformal field theory}
\author{A. V. Stoyanovsky}

\thanks{Partially supported by the RFBR grant N~04-01-00640}

\email{stoyan@mccme.ru}

\begin{abstract}
V. Drinfeld proposed conjectures on geometric Langlands correspondence
and its quantum deformation. We refine these conjectures and propose their
relationship with algebraic conformal field theory.
\end{abstract}

\maketitle

\section*{Introduction}

Geometric Langlands correspondence was proposed by V.~G.~Drinfeld [1] as a geometric analog
of the Langlands conjecture relating Galois representations of a global field with
representations of adelic algebraic group. In the geometric Langlands correspondence
the global field is replaced by the field of functions on a complete nonsingular
algebraic curve $C$ defined over the field of complex numbers $\CC$; Galois representations
are replaced by local systems on the curve $C$;
finally, representations of adelic algebraic group
are replaced by $\D$-modules (or constructible sheaves) on the moduli space of bundles on $C$.
Thus, the conjectural geometric Langlands correspondence is a correspondence between
(certain) $G^\vee$-local systems on $C$ (here $G^\vee$ is a semisimple algebraic group over
$\CC$) and between (certain) $\D$-modules on the moduli space $\Bun_G$ of principal $G$-bundles
on $C$ (here $G$ is the group Langlands dual to $G^\vee$).

The goal of this paper is to state a conjecture on deformation of the geometric Langlands
correspondence, and to relate it with constructions of algebraic conformal field theory [2].
This conjecture, originally proposed also by V.~G.~Drinfeld, looks more simple and
symmetric than the conjecture on the geometric Langlands correspondence itself. Its
relation with conformal field theory found by the author is a refinement and an argument
in favour of these conjectures, because it unifies many notions into a self-consistent picture.

For a more detailed exposition of the material of this paper, see [3].

The author is deeply grateful to B.~L.~Feigin and E.~V.~Frenkel; the main constructions of the
paper arose from discussions with them. The author is also grateful to 
A.~A.~Beilinson and V.~G.~Drinfeld for permission to cite their 
unpublished conjectures.

\section{The main conjectures}

\subsection{} Let $G$ be a simple algebraic group over $\CC$, let $C$ be a complete smooth
algebraic curve of genus $g$, and let $\Bun_G$ be the moduli space of principal $G$-bundles
on the curve $C$. It is known [3] that for any affine algebraic group $A$ the moduli space
$\Bun_A$ of principal $A$-bundles on $C$ is a smooth algebraic stack [4]. Moreover, this
stack can be covered by open substacks of the form $X_n/G_n$, where $X_n$ is a smooth
algebraic variety, and $G_n$ is an affine algebraic group acting on $X_n$. Hence various
objects (functions, $\D$-modules, sheaves, etc.) on the stack $\Bun_A$ can be
obtained by glueing corresponding $G_n$-equivariant objects on the varieties $X_n$.

Denote by $G^\vee$ the group Langlands dual to $G$. For a principal $G^\vee$-bundle
$P^\vee$ on the curve $C$, the space of (algebraic) connections on $P^\vee$ is an affine space,
whose associated vector space is the cotangent space
$$
T_{P^\vee}^*\Bun_{G^\vee}\simeq\Gamma(C,\Omega^1_C\otimes\ad P^\vee)
$$
to the stack $\Bun_{G^\vee}$ at the point $P^\vee$.
Hence we call the moduli space of $G^\vee$-local systems on $C$,
i.~e. principal $G^\vee$-bundles with a connection, by the {\it twisted cotangent bundle}
to the stack $\Bun_{G^\vee}$, and denote it by $\widetilde T^*\Bun_{G^\vee}$. It is known [3]
that the cocycle of the affine bundle $\widetilde T^*\Bun_{G^\vee}$ is obtained from the
cocycle of the canonical line bundle $\omega_{\Bun_{G^\vee}}$ on the stack $\Bun_{G^\vee}$
via the homomorphism
$$
d\log:H^1(\Bun_{G^\vee}, \OO^*_{\Bun_{G^\vee}})\to H^1(\Bun_{G^\vee}, 
\Omega^1_{\Bun_{G^\vee}}).
$$

\subsection{} {\bf Conjecture 1.} [6] The derived category of $\D_{\Bun_G}$-modules
is equivalent to the derived category of quasicoherent
$\OO_{\widetilde T^*\Bun_{G^\vee}}$-modules.

This conjecture (as well as definition of the derived categories [5]) is due to Beilinson
and Drinfeld. They refine it in the following way. The required equivalence of derived
categories should be given by the ``kernel'' $\LL_0$, an object of derived 
category of
$\D_{\Bun_G}\boxtimes\OO_{\widetilde T^*\Bun_{G^\vee}}$-modules. The 
restriction of this object
to the open substack $\widetilde T^*\Bun_{G^\vee}^\circ$ of irreducible $G^\vee$-local
systems should be a $\D_{\Bun_G}\boxtimes\OO_{\widetilde 
T^*\Bun_{G^\vee}^\circ}$-module,
flat as an $\OO_{\widetilde T^*\Bun_{G^\vee}^\circ}$-module, whose fiber 
over the local
system
$$
\P^\vee=(P^\vee,\nabla)\in\widetilde T^*\Bun_{G^\vee}^\circ
$$
should be the unique, up to isomorphism, holonomic $\D_{\Bun_G}$-module $\F_{\P^\vee}$,
which is a Hecke eigen-$\D_{\Bun_G}$-module with eigenvalue $\P^\vee$. (For a definition of
this notion, see [5].) The correspondence $\P^\vee\to\F_{\P^\vee}$ is called the geometric
Langlands correspondence.

\subsection{} Let us now proceed to the conjecture on deformation of the geometric
Langlands correspondence. Denote
$$
\xi=\omega_{\Bun_G}^{\otimes\left(-\frac1{2h^\vee}\right)}\in\Pic\Bun_G\otimes_\ZZ\QQ
$$
($\Pic$ is the Picard group), where $h^\vee$ is the dual Coxeter number of the group $G$.
In the case when $G$ is simply connected, it is known [5] that $\xi$ is the positive generator
of the group $\Pic\Bun_G\simeq\ZZ$. Similarly, introduce the element
$$
\xi^\vee=\omega_{\Bun_{G^\vee}}^{\otimes\left(-\frac1{2h}\right)}\in\Pic\Bun_{G^\vee}
\otimes_\ZZ\QQ.
$$

{\bf Conjecture 2.} [6] The derived category of twisted
$\D_{\Bun_G}(\xi^{\otimes\kappa})$-modules is equivalent to the derived category of
twisted $\D_{\Bun_{G^\vee}}(\xi^{\vee\otimes\kappa^\vee})$-modules
for any $\kappa\in\CC$, $\kappa\ne0$, where $\kappa^\vee=1/(r\kappa)$, and $r$ is the maximal
multiplicity of an edge in the Dynkin diagram of the group $G$ (or $G^\vee$);
$r=1$, $2$, or $3$.

The idea of this conjecture is due to V.~G.~Drinfeld. This conjecture also has several
refinements. We are going to state them in the rest of this paper.

The required equivalence of derived categories should be given by a kernel
$\LL_\kappa$ which is an object of derived category of
$\D_{\Bun_G}(\xi^{\otimes\kappa})
\boxtimes\D_{\Bun_{G^\vee}}(\xi^{\vee\otimes\kappa^\vee})$-modules. All refinements of
Conjecture~2 stated below deal with the properties of this kernel.

\subsection{} {\bf Property 1. Dependence on the parameter $\kappa$; classical limits.}

Let us first define the notion of asymptotic twisted $\D$-module, or
$\D_X^{\asym}(\xi^{\otimes t})$-module, on a smooth variety $X$ with a line bundle $\xi$.
This notion is introduced in [3]; let us briefly discuss it.

There exists a natural sheaf $\D_X^{\asym}(\xi^{\otimes t})$ of quasicoherent
$\OO_{\PP^1}$-algebras on the product $\PP^1\times X$, flat as an 
$\OO_{\PP^1}$-module,
whose fiber at the point $\kappa\in\PP^1$, $\kappa\ne\infty$, is isomorphic to the sheaf
$\D_X(\xi^{\otimes\kappa})$ of twisted differential operators, and the fiber over the point
$\infty\in\PP^1$ is isomorphic to the sheaf $\pi_*\OO_{\widetilde T^*X}$, 
where
$\pi:\widetilde T^*X\to X$ is the twisted cotangent affine bundle over $X$, whose
cocycle coresponds to the cocycle of the bundle $\xi$ under the homomorphism
$$
d\log:H^1(X,\OO_X^*)\to H^1(X,\Omega^1_X).
$$
Sections of the sheaf $\D_X^{\asym}(\xi^{\otimes t})$ are called twisted asymptotic
differential operators, cf. [7]. On a sufficiently small open subset $U\subset X$,
on which the bundle $\xi$ is trivial, we have
$$
\Gamma((\PP^1\setminus\{0\})\times U, \D_X^{\asym}(\xi^{\otimes t}))\simeq
\OO_U[t^{-1},t^{-1}\partial_1,\ldots,t^{-1}\partial_n],
$$
where $\partial_1$, $\ldots$, $\partial_n$ is a basis of vector fields on $U$, $t$ is the
parameter on the line $\PP^1$.

By definition, a $\D_X^{\asym}(\xi^{\otimes t})$-module is a sheaf of
$\D_X^{\asym}(\xi^{\otimes t})$-modules on the product $\PP^1\times X$ quasicoherent as an
$\OO_{\PP^1\times X}$-module.

Let us return to the kernel $\LL_\kappa$. The property of dependence of 
$\LL_\kappa$ on the
parameter $\kappa$ states that the objects $\LL_0$, $\LL_\kappa$ are the 
fibers of an object
$\LL_t$ of the derived category of
$\D_{\Bun_G}^{\asym}(\xi^{\otimes t})
\boxtimes_{\OO_{\PP^1}}\D_{\Bun_{G^\vee}}^{\asym}(\xi^{\vee\otimes 
t^\vee})$-modules
on the product $\PP^1\times\Bun_G\times\Bun_{G^\vee}$, flat as an 
$\OO_{\PP^1}$-module.
Here $t^\vee=1/(rt)$. The fiber of this object at $t=\kappa$ is the object 
$\LL_\kappa$,
the fiber at $t=0$ is the object $\LL_0$ which is the kernel of the 
geometric Langlands
correspondence defined in 1.2 above, and the fiber at $t=\infty$ is the 
object $\LL_\infty$
which is the kernel of the geometric Langlands correspondence for the group $G^\vee$, i.~e.,
the object $\LL_\infty$ is obtained from $\LL_0$ by exchanging the roles 
of the groups $G$ and
$G^\vee$.

\subsection{} {\bf Property 2: singular support of the kernel 
$\LL_\kappa$}.
One has the Hitchin map [5]
$$
\chi_G: T^*\Bun_G\to\oplus_{i=1}^{\rk G}\Gamma(C, \omega_C^{\otimes d_i}),
$$
where $d_i$ are the exponents of the group $G$, $\rk G$ is the rank of $G$. After restriction
to certain open dense subset $U\subset\Bun_G\times\Bun_{G^\vee}$, the 
kernel $\LL_\kappa$
should be a coherent $\D_{\Bun_G}(\xi^{\otimes\kappa})
\boxtimes\D_{\Bun_{G^\vee}}(\xi^{\vee\otimes\kappa^\vee})$-module whose singular support
coincides with the preimage of the diagonal
$$
\Delta_\kappa=\{v_i,v_i^\vee\in\Gamma(C,\omega_C^{\otimes d_i}):
v_i^\vee=\kappa^{d_i}v_i\}_{i=1}^{\rk G}
$$
under the product of Hitchin maps
$$
\chi_G\times\chi_{G^\vee}: T^*\Bun_G\times T^*\Bun_{G^\vee}
\to\left(\oplus_{i=1}^{\rk G}\Gamma(C, \omega_C^{\otimes d_i})\right)^{\oplus2}.
$$
The classical limit of this property as $\kappa\to0$ yields the conjecture that
the singular support of $\D_{\Bun_G}$-modules $\F_{\P^\vee}$ from the geometric Langlands
correspondence is contained in the global nilpotent cone, which is the preimage of zero
under the Hitchin map.

\section{Relation with conformal field theory}

\subsection{} For simplicity assume in this Section that the group $G$ is adjoint, and the
group $G^\vee$ is simply connected. For the general case, see [3].

Let us fix a Borel subgroup $B\subset G$ with the unipotent radical $N$; let $H=B/N$ be the
Cartan group. Consider the diagram
$$
\leqno{(*)}
\xymatrix{
       & \Bun_B^{>0}\ar[dl]_{\sigma}\ar[dr]^{\rho} &   &   &   \\
\Bun_G & & \Bun^{>0}_{B/[N,N]}\ar[dr]^{\beta} & & 
\Bun_{\omega,H}\stackrel{\iota}{\hookleftarrow}\Conf_{G^\vee}\ar[dl]_{\alpha} 
\\
       &   &   & \Bun_H^{>0} &   \\
}
$$
The only thing to be explained in this diagram is what are the spaces $\Bun_{\omega,H}$,
$\Conf_{G^\vee}$, and what does the sign $>0$ mean. To explain this, note that
$$
B/[N,N]\simeq\prod_{i=1}^{\rk G}B_i,
$$
where $B_i$ is a copy of the upper triangular Borel subgroup in $PGL(2)$. Hence
$\Bun_{B/[N,N]}$ is identified with the moduli space of sets of exact triples
$$
\leqno{(**)}\qquad\qquad\qquad\qquad 0\to\OO_C\to E_i\to L_i\to0,
$$
where $L_i$ is a line bundle on the curve $C$, $1\le i\le\rk G$. The projection
$\Bun_{B/[N,N]}\to\Bun_H$ takes a set of triples $(**)$ to the set of line bundles $(L_i)\in\Bun_H$.
By definition,
$$
\Bun_H^{>0}=\{(L_i), \deg L_i>0, 1\le i\le\rk G\};
$$
$\Bun_{B/[N,N]}^{>0}$ and $\Bun_B^{>0}$ are the preimages of $\Bun_H^{>0}$ under the natural
projections. The projection
$$
\beta:\Bun_{B/[N,N]}^{>0}\to\Bun_H^{>0}
$$
is a vector bundle whose fiber over a point $(L_i)\in\Bun_H^{>0}$ is the vector space
$$
\oplus\Ext_C^1(L_i,\OO)\simeq\oplus H^1(C,L_i^{-1}).
$$
By definition, $Bun_{\omega,H}$ is the vector bundle dual to the vector bundle $\beta$, and
$\Conf_{G^\vee}$ is an open substack of this vector bundle defined as follows:
$$
\begin{aligned}{}
\Bun_{\omega,H}&=\{(L_i,s_i),\deg L_i>0,s_i\in\Gamma(C,\omega_C\otimes L_i)\},\\
\Conf_{G^\vee}&=\{(L_i,s_i),s_i\ne0\}\simeq\{D_i: \deg D_i>2g-2\},
\end{aligned}
$$
i.~e., $\Conf_{G^\vee}$ is isomorphic to the space of sets of effective divisors $(D_i)$
on the curve $C$, i.~e. to the space of divisors with values in the semigroup $\Gamma_+$
of dominant weights of the group $G^\vee$. Let us call these sets of divisors by
{\it coloured divisors}. As a variety $\Conf_{G^\vee}$ is the disjoint union of products of
symmetric powers of the curve $C$. The fact that an open dense substack of a
vector bundle is a disjoint union of projective varieties, is due to the fact that each
point of the stack $\Bun_H$ has the group of automorphisms $(\CC^*)^{\rk G}$.

The space $\Conf_{G^\vee}$ is naturally stratified:
$$
\Conf_{G^\vee}=\sqcup_{\llambda^\vee}\Conf_{G^\vee}^{\llambda^\vee},
$$
where $\llambda^\vee=(\lambda^\vee_1,\ldots,\lambda^\vee_N)$, $\lambda^\vee_i\in\Gamma_+$, and
$$
\Conf_{G^\vee}^{\llambda^\vee}=\{\lambda_1^\vee x_1+\ldots+\lambda_N^\vee x_N, x_i\in C,
x_i\ne x_j\text{ for }i\ne j\}/\GG_{\llambda^\vee},
$$
where $\GG_{\llambda^\vee}$ is the group of permutations of indices $i$ preserving the weights
$\lambda_i^\vee$. Denote the inclusion
$\Conf_{G^\vee}^{\llambda^\vee}\hookrightarrow\Conf_{G^\vee}$ by $j_{\llambda^\vee}$.

\subsection{} {\bf Property 3 of the kernel $\LL_\kappa$.} This property 
describes the object of
derived category of twisted $\D_{\Conf_{G^\vee}^{\llambda^\vee}\times\Bun_{G^\vee}}$-modules
$$
\leqno(***)\qquad\qquad\qquad\qquad 
j_{\llambda^\vee}^!\iota^!F\rho_*\sigma^!\LL_\kappa,
$$
where all the spaces in the diagram $(*)$ are multiplied by $\Bun_{G^\vee}$; $F$ denotes the
Fourier--Laplace transform of a $\D$-module on the vector bundle $\Bun_{B/[N,N]}^{>0}$.
The object $(***)$ should be isomorphic to the twisted
$\D_{\Conf_{G^\vee}^{\llambda^\vee}\times\Bun_{G^\vee}}$-module
$\KZ_{G^\vee,\kappa^\vee}^{\llambda^\vee}$ constructed in conformal field theory.
As an
$\OO_{\Conf_{G^\vee}^{\llambda^\vee}}
\boxtimes\D_{\Bun_{G^\vee}}(\xi^{\vee\otimes\kappa^\vee})$-module,
$\KZ_{G^\vee,\kappa^\vee}^{\llambda^\vee}$ coincides with the induced module from the
$\OO_{\Conf_{G^\vee}^{\llambda^\vee}\times\Bun_{G^\vee}}$-module whose 
fiber at the point
$$
(\lambda^\vee_1x_1+\ldots+\lambda_N^\vee x_N, P^\vee)
\in\Conf_{G^\vee}^{\llambda^\vee}\times\Bun_{G^\vee}
$$
is the tensor product
$(V_{P^\vee}^{\lambda^\vee_1})_{x_1}\otimes\ldots\otimes (V_{P^\vee}^{\lambda^\vee_N})_{x_N}$,
where $V_{P^\vee}^{\lambda^\vee_i}$ is the vector bundle on $C$ associated with the
principal $G^\vee$-bundle $P^\vee$ and with the $G^\vee$-module $V^{\lambda^\vee_i}$
with highest weight $\lambda_i^\vee$. Further, this
$\OO_{\Conf_{G^\vee}^{\llambda^\vee}}
\boxtimes\D_{\Bun_{G^\vee}}(\xi^{\vee\otimes\kappa^\vee})$-module has a natural structure
of a twisted $\D_{\Conf_{G^\vee}^{\llambda^\vee}}$-module [2], which is a direct
generalization of the Knizh\-nik--Za\-mo\-lod\-chi\-kov connection to curves of genus $g$.
A direct check with the use of Riemann-Roch theorem shows that the object $(***)$ has the same
twist. The property~3 of the kernel $\LL_\kappa$ states that these two 
objects are isomorphic.

\subsection{} {\bf Classical limits of the property 3.}

a) {\it The limit $\kappa\to0$} amounts to the geometric analog of the
Casselman--Shalika--Shintani formula for the Whittaker function of an automorphic form. This
statement is that the $\D_{\Conf_{G^\vee}^{\llambda^\vee}}$-module
$$
j_{\llambda^\vee}^!\iota^!F\rho_*\sigma^!\F_{\P^\vee}
$$
is isomorphic to the local system whose fiber over the point
$\lambda^\vee_1x_1+\ldots+\lambda^\vee_Nx_N\in\Conf_{G^\vee}^{\llambda^\vee}$ equals
$(V_{\P^\vee}^{\lambda^\vee_1})_{x_1}\otimes\ldots\otimes (V_{\P^\vee}^{\lambda^\vee_N})_{x_N}$.

b) {\it The limit $\kappa\to\infty$} amounts to the construction of the
$\D_{\Bun_{G^\vee}}$-module $\F_\P$, for $\P\in\widetilde T^*\Bun_G^\circ$,
by means of a $G$-oper with regular singularities with trivial monodromy [1].

For details, see [3].

\subsection{} The following natural question arises: what is the result of applying the same
operation $j_{\llambda^\vee}^!\iota^!F\rho_*\sigma^!$ to the twisted $\D$-module
$\KZ_{G,\kappa}^\llambda$ on the product $\Conf_G^\llambda\times\Bun_G$? The arising twisted
$\D$-module on the product $\Conf_G^\llambda\times\Conf_{G^\vee}^{\llambda^\vee}$ should
be related with the $W$-algebra $W_G^\kappa\simeq W_{G^\vee}^{\kappa^\vee}$ [8]. In the case
$G=SL(2)$ it is the Virasoro algebra.

For a statement of this kind for a curve $C$ of genus $g=0$ with marked points and for
$G=SL(2)$, see [9].

\end{document}